\documentclass{amsart}

\usepackage{amssymb}
\usepackage{amsmath}
\usepackage{amsfonts}
\usepackage{latexsym}
\usepackage{amsthm}
\usepackage{amscd}
\usepackage{enumerate}
\usepackage[latin1]{inputenc}

\newtheorem{thm}{Theorem}[section]
\newtheorem{cor}[thm]{Corollary}
\newtheorem{lem}[thm]{Lemma}
\newtheorem{prop}[thm]{Proposition}
\theoremstyle{definition}
\newtheorem{defn}[thm]{Definition}
\theoremstyle{remark}

\newtheorem{rem}[thm]{Remark}
\newtheorem{exam}[thm]{Example}

 \DeclareMathOperator{\RE}{Re}
 \DeclareMathOperator{\IM}{Im}

 \newcommand{\rank}{\ensuremath{\mathrm{rank}}}
 \newcommand{\R}{\ensuremath{\mathbb{R}}}
 \newcommand{\C}{\ensuremath{\mathbb{C}}}

 \newcommand{\derp}[2]{\ensuremath{\frac{\partial #1}{\partial #2}}}
 \newcommand{\derpp}[1]{\ensuremath{\derp{}{#1}}}
 \newcommand{\sep}{\setlength{\itemsep}{0pt}}
\begin{document}

\title[Complex and CR-structures on compact Lie groups]
{Complex and CR-structures on compact Lie groups associated to Abelian
actions}

\author{J.-J. Loeb, M. Manjar\'{\i}n, M. Nicolau}

\address{J.-J. Loeb. Departament de Math\'{e}matiques et Informatique,  Universit\'e d'Angers, 2 Boulevard Lavoisier, Angers 49045,
FRANCE} \address{M. Manjar\'{\i}n. Departament d'\`{A}lgebra i Geometria,
Universitat de Barcelona, Gran Via de les Corts Catalanes 585,
Barcelona 08007, SPAIN} \address{M. Nicolau. Departament de
Matem\`{a}tiques, Universitat Aut\`{o}noma de Barcelona, Bellaterra 08193,
SPAIN}

\medskip

\email{Jean-Jacques.Loeb@univ-angers.fr, manjarin@ub.edu,
nicolau@mat.uab.es}

\thanks{This work was partially supported by the Ministerio de Educaci\' on y Ciencia,
grants MTM2004-00566 and PR2006-0038, of Spain.}




\dedicatory{}

\commby{}


\begin{abstract}
It was shown by Samelson \cite{Samel} and Wang \cite{Wang2}
that each compact Lie group $\mathrm{K}$ of even dimension
admits left-invariant complex structures. When $\mathrm{K}$
has odd dimension it admits a left-invariant CR-structure
of maximal dimension. This has been proved recently by
Charbonnel and Khalgui \cite{Charb} who have also given a complete
algebraic description of these structures.

In this article we present an alternative and more geometric construction
of this type of invariant structures on a compact Lie group
$\mathrm{K}$ when it is semisimple.
We prove that each left-invariant complex structure, or
each CR-structure of maximal dimension with a transverse CR-action
by $\mathbb R$, is induced
by a holomorphic ${\mathbb C}^l$-action on a quasi-projective
manifold $\mathrm X$ naturally associated to $\mathrm{K}$.
We then show that $\mathrm X$ admits more general
Abelian actions, also inducing complex or CR-structures
on $\mathrm{K}$  which are generically non-invariant.
\end{abstract}

\maketitle

\section{Introduction}

Each compact Lie group $\mathrm{K}$ of even dimension can be
endowed with a left-invariant complex structure. This was proved
independently by Samelson \cite{Samel} and Wang \cite{Wang2}.
At that time their construction provided one of the first known families
of compact non-K\"{a}hler complex manifolds since, by topological
reasons, a compact complex Lie group cannot be K\"{a}hler unless
it is Abelian.

Recently Charbonnel and Khalgui have given an algebraic description
of all left-invariant CR-structures of maximal dimension on compact
Lie groups  (cf. \cite{Charb}). If $\dim \mathrm{K}$ is even these structures
are invariant complex structures on $\mathrm{K}$ and there is only one
type of them. In the odd-dimensional case CR-structures of
maximal dimension are divided into two types. This classification is made
in terms of Cartan subalgebras and root systems.

The aim of this paper is to present an alternative and more geometrical
construction of left-invariant complex structures and of CR-structures
of maximal dimension admitting a transverse CR-action of $\mathbb R$
(the precise definition is given in \S 2)
on a compact Lie group $\mathrm{K}$. If $\dim \mathrm{K}$ is odd
these kind of CR-structures corresponds precisely to one of the two
types considered by Charbonnel and Khalgui.

Assume that $\mathrm{K}$ is semisimple and let $\mathrm{G}$ denote its
universal complexification. Then $\mathrm{K}$ is naturally
embedded into the quasi-projective manifold $\mathrm{G/U}$ where
$\mathrm{U}$ is a maximal unipotent complex subgroup of $\mathrm{G}$.
If $\dim \mathrm{K}=2n$, a left-invariant
complex structure on $\mathrm{K}$ is determined by a complex subgroup
$\mathrm{L}$ of $\mathrm{G}$ of dimension $n$ and such that
$\mathrm{L}\cap\mathrm{K} = \{e\}$. It turns out  that $\mathrm{L}$
is closed and solvable and, up to conjugation, it can be chosen in such
a way that $\mathrm{U}\subset\mathrm{L}$. The choice of $\mathrm{L}$
determines a ${\mathbb C}^l$-action on $\mathrm{G/U}$,
where $l=\rank\,\mathrm{K}/2$, which is transverse to $\mathrm{K}$
and therefore defines a complex structure on the compact Lie group.
In fact $\mathrm{K}$ can be identified to the orbit space of the
${\mathbb C}^l$-action. It is proved that each left-invariant complex
structure is obtained in this way.  This procedure is similar
to the construction of non-K\"{a}hler manifolds obtained as orbit spaces
of holomorphic group actions developed by Loeb-Nicolau in
\cite{MNLoeb1}, L\'{o}pez de Medrano-Verjovsky in \cite{LMVer} and
Meersseman in \cite{Meer}.

If $\dim \mathrm{K}=2n+1$, a left-invariant CR-structure of maximal dimension
is determined also by a complex subgroup $\mathrm{L}$ of $\mathrm{G}$
of dimension $n$ with $\mathrm{L}\cap\mathrm{K} = \{e\}$. However, in this case
$\mathrm{L}$ does not necessarily contain a maximal unipotent subgroup
of $\mathrm{G}$. In fact $\mathrm{L}$ contains such a subgroup $\mathrm{U}$
precisely when the CR-structure on $\mathrm{K}$ admits an invariant transverse
CR-action of $\mathbb R$. This is equivalent to the existence of a
$n+1$-dimensional complex subgroup $\mathrm{L}'$ of $\mathrm{G}$
containing $\mathrm{L}$ as a normal subgroup. Then the transverse
CR-action is given by the one-parameter subgroup of $\mathrm{K}$ defined
as $\mathrm{K}\cap\mathrm{L}'$. In this situation the construction above
can be adapted to the odd-dimensional case and we are able to describe
all the invariant CR-structures
of maximal dimension with a transverse $\mathbb R$-action on a
semisimple Lie group $\mathrm{K}$  in terms of appropriate Abelian actions
on $\mathrm{G/U}$.
The existence of invariant structures of these types on an arbitrary compact
Lie group follows easily from the above constructions.

The main interest of the geometric construction that we present here is that,
by deformation of the action, one can define structures of
the above types on $\mathrm{K}$ which are not
invariant.
In section~4 we introduce a larger class of ${\mathbb C}^l$-actions
on the quotient space $\mathrm{G/U}$ showing that they define complex structures
on $\mathrm{K}$, or CR-structures with transverse $\mathbb R$-action,
according to the parity of the dimension of the group. We also show
that the structures so obtained are generically non-invariant.


\section{Invariant complex and CR-structures on Lie groups}\label{sec:consinv}

Throughout all the paper K will be a compact connected real Lie group and
$\mathfrak{k}$ will denote its Lie algebra. Recall that the dimension and the rank
of $\mathrm{K}$ have the same parity.
We set $\rank \,
\mathrm{K}=2r$ when $\dim\mathrm{K}=2n$ and $\rank \,
\mathrm{K}=2r+1$ when $\dim\mathrm{K}=2n+1$. Let $\mathrm{G}$
denote the universal complexification of the Lie group $\mathrm{K}$
and $\mathfrak{g}=\mathfrak{k}^{\C}:=\mathfrak{k}\otimes \C$ its Lie
algebra. $\mathrm G$ is a linear algebraic group.

\begin{prop}\label{prop:liecr}
An (integrable) left-invariant CR-structure over $\mathrm{K}$ is
defined by a complex subalgebra $\mathfrak{l}$ of $\mathfrak{g}$
such that $\mathfrak{l}\cap \mathfrak{k}=\{0\}$.
\end{prop}

\begin{proof}
Notice that a left-invariant CR-structure on $\mathrm{K}$ is
determined by a complex subspace $\mathfrak{l}$ of $\mathfrak{g}$,
which defines the vectors of $T_{e}^{\C}\mathrm{K}$ of type $(0,1)$, fulfilling
the condition $\mathfrak{l}\cap \overline{\mathfrak{l}}=\{0\}$. This condition is
equivalent to $\mathfrak{l}\cap \mathfrak{k}=\{0\}$. On the other hand the
CR-structure defined by $\mathfrak{l}$ is integrable (or involutive)
if and only if $[\mathfrak{l},\mathfrak{l}]\subset \mathfrak{l}$,
i.e., if $\mathfrak{l}$ is a subalgebra.
\end{proof}

In this article we will be concerned with (left-invariant) CR-structures
of {\sl maximal dimension} on $\mathrm{K}$; that is, $\dim_{\C}\mathfrak{l}=n$ if
$\dim_{\R}\mathrm{K}=2n$ or $\dim_{\R}\mathrm{K}=2n+1$.
Notice that in the even-dimensional case, i.e. $\dim_{\R}\mathrm{K}=2n$, a
CR-structure
of maximal dimension is nothing but a complex structure on $\mathrm{K}$.
By convention the subalgebra
$\mathfrak{l}$ will always correspond to the distribution of
vector fields of type (0,1) of the CR-structure.

By a CR-structure with a {\sl transverse $\R$-action} on an odd-dimensional manifold $\mathrm{M}^{2n+1}$ we mean a CR-structure of maximal
dimension  on $\mathrm{M}$ defined by a complex subbundle
$\Phi^{0,1}\subset T^{\C}\mathrm{M}$ together with a CR-action of $\R$
induced by a vector field $\xi$ which is transverse to
$\mathcal D = T\mathrm{M}\cap (\Phi^{0,1}\oplus\overline{\Phi^{0,1}})$
at each point, i.e. $T\mathrm{M} = \langle\xi\rangle\oplus \mathcal D$.

The notion of CR-structure with a transverse $\R$-action on an odd-dimensional
manifold
coincides with that of \emph{normal almost contact structure} (cf. \cite{Blair}).
Such a structure will be called \emph{nacs} for shortness. The following proposition
gives a characterization of left-invariant nacs on $\mathrm{K}$.

\begin{prop}
A left-invariant nacs over a compact connected real Lie group
$\mathrm{K}$ such that $\dim_{\R}\mathrm{K}=2n+1$ is determined
by a pair of complex subalgebras $\mathfrak{l}\subset \mathfrak{l}'$
of $\mathfrak{g}$ of complex dimension $n$ and $n+1$ respectively
such that:
\begin{enumerate}[\bf (a)] \sep
\item $\mathfrak{l}\cap \mathfrak{k}=\{0\}$;
\item $\mathrm{dim}_{\R}\mathfrak{l}'\cap \mathfrak{k}=1;$
\item $\mathfrak{l}$ is an ideal of $\mathfrak{l}'$, i.e. $[\mathfrak{l},\mathfrak{l}']\subset \mathfrak{l}$.
\end{enumerate}
\end{prop}

\begin{proof}
Proposition \ref{prop:liecr} says that the subalgebra
$\mathfrak{l}$ defines a left-invariant CR-struc\-tu\-re on
$\mathrm{K}$. Note that $\mathfrak{l}'\cap \mathfrak{k}=\langle \xi
\rangle_{\R}$ corresponds to the left-invariant vector field
defining the CR-action. Clearly the vector field $\xi$ is transverse
to the CR-structure determined by $\mathfrak{l}$ and condition \textbf{(c)}
implies that it induces a CR-action.

Conversely, if a normal almost contact structure on $\mathrm{K}$
is given by a subalgebra $\mathfrak{l}$ defining the CR-structure
and a left-invariant vector field $\xi$ defining the CR-action
then it is enough to set
$\mathfrak{l}'=\mathfrak{l}\oplus \langle \xi \rangle_{\C}$. Notice
that $\mathfrak{g}=\mathfrak{l}\oplus \overline{\mathfrak{l}}\oplus
\langle \xi \rangle_{\C}$ and therefore
$\mathrm{dim}_{\R}\mathfrak{l}'\cap \mathfrak{k}=1$.
\end{proof}

\begin{rem}\label{rem:nacsxs1}
If the complex subalgebras $\mathfrak{l}\subset
\mathfrak{l}'=\mathfrak{l}\oplus \langle \xi \rangle_{\C}$ of $\mathfrak{g}$
determine a left-invariant nacs over $\mathrm{K}$ then the product
$\mathrm{K}\times S^1$
can be endowed with a left-invariant complex structure in a natural way.
It suffices to define the distribution of vector fields of type
$(0,1)$ on $\mathrm{K}\times S^1$ as the involutive distribution
generated by the subalgebra $\mathfrak{l}\oplus \langle
\xi+i\derpp{t}\rangle$, where $\derpp{t}$ is the vector field determined
by the natural $S^1$-action.

In a similar way if ${\mathrm K}_1$ and ${\mathrm K}_2$
are compact connected Lie groups with left-invariant nacs defined by
pairs of subalgebras $\mathfrak{l}_1\subset
\mathfrak{l}'_1=\mathfrak{l}_1\oplus \langle \xi_1 \rangle_{\C}$
and $\mathfrak{l}_2\subset
\mathfrak{l}'_2=\mathfrak{l}_2\oplus \langle \xi_2 \rangle_{\C}$
respectively then the product ${\mathrm K}_1 \times {\mathrm K}_2$
is endowed with a left-invariant
complex structure defined by the subalgebra
$\mathfrak{l}_1\oplus\mathfrak{l}_2 \oplus \langle
\xi_1+i\xi_2\rangle$ of $\mathfrak{g}_1\oplus\mathfrak{g}_2$.
\end{rem}

Charbonnel and Khalgui, in a recent paper \cite{Charb}, have given a
complete characterization of left-invariant CR-struc\-tu\-res of
maximal dimension on a compact Lie group $\mathrm{K}$ in terms of
the algebraic structure of the subalgebra $\mathfrak{l}$. According
to this structure they divide these CR-structures into two types:
CR0 and CRI. With this terminology left-invariant nacs on
$\mathrm{K}$ correspond exactly to CR-structures of type CR0.

In the characterization given by Charbonnel and Khalgui a crucial
fact is the solvability of the subalgebra  $\mathfrak{l}$. Using a
more geometrical approach we give here an independent proof of this
result for left-invariant complex structures and for nacs when the
Lie group $\mathrm{K}$ is semisimple. The interest of this
alternative point of view is two-folded. On one hand we are able to
prove that, under the above hypothesis, the connected Lie subgroup
$\mathrm{L}$ of $\mathrm{G}$ associated to the subalgebra
$\mathfrak{l}$ is closed. On the other hand, this approach allows us
to construct, in the last section, non-invariant complex structures
and nacs on $\mathrm{K}$.

\begin{rem}\label{rem:generalcase}
Recall that a compact Lie group $\mathrm{K}$ has a finite covering group
which is the product of a semisimple compact Lie group $\mathrm{K}_1$
with a torus $(S^1)^m$. Thus one can write
$$
\mathrm{K} = \Gamma\backslash(\mathrm{K}_1\times(S^1)^m)
$$
where $\Gamma$ is a finite subgroup of the center $Z(\mathrm{K}_1\times(S^1)^m) =
Z(\mathrm{K}_1)\times(S^1)^m$ of $\mathrm{K}_1\times(S^1)^m$. By this reason we
will start discussing the semisimple case and at the end we will prove the
existence of complex structures or nacs on an arbitrary compact Lie group
using the above description.
\end{rem}

\begin{thm}\label{teo:principal}
Let $\mathrm{K}$ be a compact connected semisimple Lie group
endowed with a left-invariant
complex structure defined by a subalgebra $\mathfrak{l}$ of
$\mathfrak{g}$ or with a left-invariant nacs defined by a pair of subalgebras
$\mathfrak{l}\subset \mathfrak{l}'$ of $\mathfrak{g}$. Let $\mathrm{L}$ be the connected
complex Lie subgroup associated to $\mathfrak{l}$. Then the subalgebra $\mathfrak{l}$
is solvable and the subgroup $\mathrm{L}$ is closed.

Moreover, the inclusion $\mathrm{K}\subset \mathrm{G}$ induces a CR-embedding
$\mathrm{K}\hookrightarrow \mathrm{G}/\mathrm{L}$ of $\mathrm{K}$ into the homogeneous
complex manifold $\mathrm{G}/\mathrm{L}$ and
\begin{enumerate}[\bf (i)]\sep
\item if $\dim_{\R}\mathrm{K}=2n$, the map $\mathrm{K}\hookrightarrow \mathrm{G}/\mathrm{L}$
is a biholomorphism,
\item if $\dim_{\R}\mathrm{K}=2n+1$ then $\mathrm{G}/\mathrm{L}$ is
diffeomorphic to the product $\mathrm{K}\times\mathbb R$ and
the map $\mathrm{K}\hookrightarrow \mathrm{G}/\mathrm{L}$ identifies $\mathrm{K}$ with
$\mathrm{K}\times \{0\}$. Furthermore there is a holomorphic vector field $\zeta$ on
$\mathrm{G}/\mathrm{L}$ which is transverse to $\mathrm{K}$ and whose real part
is tangent to $\mathrm{K}$ and defines the CR-action of the nacs on $\mathrm{K}$.
\end{enumerate}
\end{thm}

\begin{cor}\label{cor:principal}
Let $\mathrm{K}$ be a compact connected semisimple Lie group
endowed with a nacs defined by a
pair of subalgebras $\mathfrak{l}\subset \mathfrak{l}'$ of
$\mathfrak{g}$. Then $\mathfrak{l}'$ is solvable.
\end{cor}

\begin{proof}(of the Corollary)
Since $\mathfrak{l}'=\mathfrak{l}\oplus \langle \xi \rangle_{\C}$
and $\mathfrak{l}$ is an ideal of $\mathfrak{l}'$ we have
$[\mathfrak{l}',\mathfrak{l}']\subset \mathfrak{l}$, thus
$\mathfrak{l}'$ is solvable if and only if $\mathfrak{l}$ is
solvable.
\end{proof}

\begin{proof}(of the Theorem)
We consider the compact Lie group $\widehat{\mathrm{K}}$ endowed with
a left-invariant complex structure defined as follows. If  $\dim_{\R}\mathrm{K}=2n$
we set $\widehat{\mathrm{K}} = \mathrm{K}$ and if $\dim_{\R}\mathrm{K}=2n+1$
we define $\widehat{\mathrm{K}}$ as the product $\mathrm{K}\times S^1$ with
the left-invariant complex structure defined in remark \ref{rem:nacsxs1}.

We shall first prove that there exist a complex closed subgroup
$\mathrm{L}$ of the universal complexification $\mathrm{G}$ of
$\mathrm{K}$ and an embedding K$\hookrightarrow$G/L
such that the complex or CR-structure on K defined by
$\mathfrak{l}$ agrees with the one induced by the embedding and the
homogeneous complex structure on G/L. Next we will see that the
Lie subalgebra associated to $\mathrm{L}$ coincides with $\mathfrak{l}$.

We consider the complex Lie group
$\widehat{\mathrm{G}}=\mathrm{Aut}_{\C}(\widehat{\mathrm{K}})$
of holomorphic automorphisms of the complex manifold $\widehat{\mathrm{K}}$
and we denote by $\widehat{\mathrm{L}}$ the isotropy group of
$e\in \widehat{\mathrm{K}}$, i.e. $\widehat{\mathrm{L}}=
\{f\in \widehat{\mathrm{G}}\, : \, f(e)=e\}$. Then $\widehat{\mathrm{L}}$ is closed
and, since $\widehat{\mathrm{K}}$ acts holomorphically and transitively onto itself by
left-translations, there is an inclusion  $\widehat{\mathrm{K}}\subset \widehat{\mathrm{G}}$
and the complex manifold $\widehat{\mathrm{K}}$ is naturally identified to
$\widehat{\mathrm{G}}/\widehat{\mathrm{L}}$.
In particular there is a CR-embedding
$\mathrm{K}\hookrightarrow \widehat{\mathrm{G}}/\widehat{\mathrm{L}}$
which is just a biholomorphism in the even-dimensional case.

Now we will see that $\mathfrak{k}$  is totally real
in the complex Lie algebra $\widehat{\mathfrak{g}}$ of $\widehat{\mathrm{G}}$.
Let us consider the ideal $\mathfrak{m}$ of $\mathfrak{k}$ defined by
$\mathfrak{m}=\mathfrak{k} \cap \mathrm{i}\mathfrak{k}$, where i
denotes the product by $\sqrt{-1}$ inside the
complex Lie algebra $\widehat{\mathfrak{g}}$.
Since K is
semisimple we have $\mathfrak{k}\cong \mathfrak{i}_1 \oplus ...
\oplus \mathfrak{i}_p$ where $\mathfrak{i}_1,...,\mathfrak{i}_p$
are simple ideals of $\mathfrak{k}$ and we can write
$\mathfrak{m}$ as a sum of some of these ideals. In particular
$[\mathfrak{m},\mathfrak{m}]=\mathfrak{m}$. On the other hand
$\mathfrak{m}$ is a complex subalgebra and the associated Lie
subgroup $\mathrm{M}$ of $\mathrm{K}$ is a complex Lie group.
Since $\mathfrak{m}$ is an ideal of $\mathfrak{k}$ then $\mathrm{M}$
is a normal subgroup of the semisimple compact Lie group K
and therefore it is closed. It follows that $\mathrm{M}$ is a compact complex Lie group
and consequently Abelian.
Thus $\mathfrak{m} = [\mathfrak{m},\mathfrak{m}] = 0$ and $\mathfrak{k}$ is
totally real in $\widehat{\mathfrak{g}}$.

Let $\mathrm{G}'$ be the connected complex subgroup of $\widehat{\mathrm{G}}$
associated to the complex Lie subalgebra
$\mathfrak{g}:=\mathfrak{k}\oplus \mathrm{i}\, \mathfrak{k}$ of
$\widehat{\mathfrak{g}}$. Since $\mathrm{K}$ is totally real in the
complex Lie group $\mathrm{G}'$ and the
Lie algebra of $\mathrm{G}'$ is a complexification of the Lie
algebra of $\mathrm{K}$ then $\mathrm{K}$ is a compact real form of $\mathrm{G}'$.
Since $\mathrm{G}'$ is semisimple $\mathrm{K}$ is a maximal compact subgroup
of $\mathrm{G}'$. On the other hand $\mathrm{K}$ is also a maximal compact subgroup
of its universal complexification $\mathrm{G}$.
By the universal property of $\mathrm{G}$ there is a group morphism
$\mathrm{G}\rightarrow \mathrm{G}'$ sending $\mathfrak{k}$ identically onto itself.
This morphism is therefore a covering map and since $\mathrm{K}$ is
as a deformation retract for both $\mathrm{G}$ and $\mathrm{G}'$ one deduces that $\mathrm{G}\cong\mathrm{G}'$.

Let $\mathrm{L}_{1}$ denote the closed complex subgroup
$\widehat{\mathrm{L}}\cap \mathrm{G}$ of $\mathrm {G}$. The quotient
$\mathrm{G}/\mathrm{L}_{1}$ is naturally included into
$\widehat{\mathrm{G}}/\widehat{\mathrm{L}}$ and there is a CR-embedding
\begin{equation}\label{embed}
\mathrm{K} \hookrightarrow \mathrm{G}/\mathrm{L}_{1} \subset
\widehat{\mathrm{G}}/\widehat{\mathrm{L}}\cong \widehat{\mathrm{K}}
\end{equation}
Let us see that in fact $\mathrm{G}/\mathrm{L}_{1} =
\widehat{\mathrm{G}}/\widehat{\mathrm{L}}\cong \widehat{\mathrm{K}}$.
By construction the composition of all the inclusions in (\ref{embed}) is just the identity of
$\mathrm{K} = \widehat{\mathrm{K}}$ when $\dim_{\R}\mathrm{K}= 2n$
and identifies $\mathrm{K}$ with $\mathrm{K}\times\{e\}\subset\widehat{\mathrm{K}}$
when  $\dim_{\R}\mathrm{K} = 2n +1$. The assertion is thus clear in the
even-dimensional case. Since $\mathrm{G}/\mathrm{L}_{1} $ is a complex manifold
we also deduce from (\ref{embed}) that
$\dim_{\C}\mathrm{G}/\mathrm{L}_{1}=n+1$ if $\dim_{\R}\mathrm{K} = 2n +1$. In this
last case the map
\begin{equation}\label{fibration}
\psi: \mathrm{G} \rightarrow \widehat{\mathrm{K}},
\end{equation}
obtained as the composition
$\psi: \mathrm{G} \rightarrow  \mathrm{G}/\mathrm{L}_{1} \subset
\widehat{\mathrm{G}}/\widehat{\mathrm{L}}\cong \widehat{\mathrm{K}}$,
is surjective showing the desired identity. In fact, since $\widehat{\mathrm{K}}$
is an exponential Lie group, any element $h\in \widehat{\mathrm{K}}$ can be written
$h = \exp(v)$, with $v$ a left-invariant vector field, and then $h = \psi(\exp(\tilde{v}))$
where $\tilde{v}$ is a vector fulfilling $d_0{\psi}(\tilde{v}) = v$.

We deduce that, in both cases, the map (\ref{fibration}) defines a fibration with fibre $\mathrm{L}_{1}$.
The corresponding homotopy long exact sequence gives
\begin{equation}\label{homotopy}
0 = \pi_{2}(\widehat{\mathrm{K}}) \rightarrow \pi_{1}(\mathrm{L}_{1}) \rightarrow
\pi_{1}(\mathrm{G})\overset{\varphi}\rightarrow \pi_{1}(\widehat{\mathrm{K}}) \rightarrow
\pi_{0}(\mathrm{L}_{1}) \rightarrow 0.
\end{equation}
If $\dim_{\R}\mathrm{K} = 2n$ the map $\varphi: \pi_{1}(\mathrm{G})\rightarrow
\pi_{1}(\widehat{\mathrm{K}})$ is an isomorphism because $\widehat{\mathrm{K}} =
\mathrm{K}$ is a deformation retract of $\mathrm{G}$ and we deduce that
$\mathrm{L}_{1}$ is connected. In this case we define $\mathrm{L} = \mathrm{L}_{1}$.
If $\dim_{\R}\mathrm{K} = 2n+1$, since $\psi(\mathrm{K}) = \mathrm{K}\times\{e\}$
the morphism $\varphi$ is
injective with a cokernel isomorphic to $\pi_{0}(\mathrm{L}_{1})\cong\mathbb Z$.
In this case we define $\mathrm{L}$ as the connected component of
the identity of $\mathrm{L}_{1}$. Notice that $\mathrm{L}$ is closed in $\mathrm{G}$
and that $\mathrm{G}/\mathrm{L}$ is
diffeomorphic to the product $\mathrm{K}\times\mathbb R$. Moreover the
inclusion $\mathrm{K} \subset \mathrm{G}/\mathrm{L}_{1}$ induces a
well defined CR-embedding $\mathrm{K} \hookrightarrow \mathrm{G}/\mathrm{L}$.

Therefore $\mathrm{K}$ is biholomorphic to $\mathrm{G}/\mathrm{L}$
in the even-dimensional case and, in the odd-dimensional case, $\mathrm{K}$
inherits its CR-structure
from the complex structure of $\mathrm{G}/\mathrm{L}$ by means of the embedding
$\mathrm{K}\hookrightarrow \mathrm{G}/\mathrm{L}$.
In both cases this is equivalent to say that
the Lie algebra of L coincides with $\mathfrak{l}$.
This follows from the identities
$\mathrm{TG} = \mathfrak{k}\oplus\mathrm{i}\mathfrak{k} = \mathfrak{l} \oplus \bar{\mathfrak{l}}$
and $\mathrm{TG} = \mathfrak{k}\oplus\mathrm{i}\mathfrak{k} =
\mathfrak{l} \oplus \bar{\mathfrak{l}}\oplus \langle\xi\rangle$ in the even-dimensional case
and in the odd-dimensional case respectively.

If $\dim_{\R}\mathrm{K} = 2n +1$ let $\mathrm{L}'$ be the connected
subgroup associated to $\mathfrak{l}' = \mathfrak{l}\oplus\langle\xi\rangle$.
Since $\mathrm{L}$
is normal in $\mathrm{L}'$ the element $\xi\in \mathfrak{l}'$ induces a well
defined vector field $\zeta$ on  $\mathrm{G}/\mathrm{L}$ with the properties
stated in the theorem.

Finally we end the proof by showing the solvability of the algebra
$\mathfrak{l}$. By Levi-Malcev's theorem the complex algebra $\mathfrak{l}$
can be decomposed
as a sum of a solvable ideal $\mathfrak{r}$ and a semisimple
subalgebra $\mathfrak{s}$. If $\mathfrak{l}$ is not solvable then
the complex subalgebra $\mathfrak{s}\neq 0$ admits a real compact form T,
i.e. there exists a compact real Lie subgroup T of L with Lie
algebra $\mathfrak{t}$ such that $\mathfrak{s}=\mathfrak{t}^{\C}$.
As K is a maximal compact subgroup of G there exists an element
$g\in \mathrm{G}$ such that $\mathrm{T}\subset g \mathrm{K}
g^{-1}=\mathrm{K}'$, in particular K$'\cap \mathrm{L}\supseteq
\mathrm{T}\neq \{e\}$. Since the action of K over
$\mathrm{G}/\mathrm{L}$ by left-translations is free we deduce
that the action of K$'$ by left-translations over
$\mathrm{G}/\mathrm{L}$ is also free. It
follows that K$'\cap$L$=\{e\}$. This leads to a
contradiction and therefore $\mathfrak{l}$ must be solvable.
\end{proof}

\begin{rem} Given a complex subalgebra $\mathfrak{l}$ with
$\mathrm{dim}_{\C}\mathfrak{l}=n$ and
$\mathfrak{l}\cap \mathfrak{k}=\{0\}$ it does not always exist a
subalgebra $\mathfrak{l}'\supset \mathfrak{l}$ defining a nacs.
The Lie groups $\mathrm{SO(3)}$ and
$\mathrm{SU(2)}$ provide examples of this situation, as we will
see in the next section.
\end{rem}

Since the algebras $\mathfrak{l}$ and $\mathfrak{l}'$ are solvable
there is a Borel subgroup $\mathrm B$ of $\mathrm G$, that is a maximal
solvable algebraic subgroup, whose Lie
algebra  $\mathfrak{b}$ contains $\mathfrak{l}$ and $\mathfrak{l}'$.
The next statement describes the relation between $\mathfrak{l}$,
$\mathfrak{l}'$ and $\mathfrak{b}$.

\begin{thm}\label{teo:borel}
Let $\mathrm{K}$ be a compact connected semisimple Lie group with
Lie algebra $\mathfrak{k}$ endowed with a left-invariant complex
structure defined by a subalgebra $\mathfrak{l}$ of $\mathfrak{g}$
or with a nacs defined by a pair $\mathfrak{l}\subset \mathfrak{l}'$ of
subalgebras of $\mathfrak{g}$. Let $\mathrm B$ be a Borel subgroup of
$\mathrm G$ whose Lie algebra $\mathfrak{b}$ contains  $\mathfrak{l}$,
and $\mathfrak{l}'$ in the odd-dimensional case,
and set $\mathfrak{u}=[\mathfrak{b},\mathfrak{b}]$. Then there exist
a Cartan subalgebra
$\mathfrak{r}=\mathfrak{t}^{\C}$ contained in $\mathfrak{b}$ such that
\begin{enumerate}[\bf (i)]\sep
\item $\mathfrak{t}$ is a maximal Abelian subalgebra of $\mathfrak{k}$,
\item $\mathfrak{l}=(\mathfrak{l}\cap \mathfrak{r})\oplus \mathfrak{u}$,
and also $\mathfrak{l}'=(\mathfrak{l}'\cap \mathfrak{r})\oplus
\mathfrak{u}$ in the case of nacs. In particular $\mathfrak{u}\subset \mathfrak{l}$.
\end{enumerate}
\end{thm}

We begin by stating some results that will be used in the proof of
the theorem.

\begin{lem}\label{lem:u2}
Every connected complex subgroup \emph{M} of $(\C^{\ast})^k$ is
isomorphic to $\C^l\times (\C^{\ast})^s$ for $l,s\geq 0$ and
$l+s\leq k$. Moreover:
\begin{enumerate}[\bf (a)]\sep
\item If $\dim_{\R}(\mathrm{M}\cap(S^1)^{k})=0$ then \emph{M} is
isomorphic to $\C^l$ and $0\leq l\leq k/2$.
\item If $\dim_{\R}(\mathrm{M}\cap(S^1)^{k})=1$ then \emph{M} is
isomorphic to $\C^l\times\C^{\ast}$ or $\C^{l+1}$ and $0\leq l
\leq (k-1)/2$.
\end{enumerate}
\end{lem}

\begin{proof}
By a theorem of Morimoto (cf. \cite{Morimo2}) a Lie subgroup
$\mathrm{M}$ of $(\C^{\ast})^k$ is isomorphic to
$\mathrm{M}^0\times\C^l \times (\C^{\ast})^s$ where $\mathrm{M}^0$
is a (HC)-group, that is an Abelian group without holomorphic functions other
than the constants. As $(\C^{\ast})^k$ is Stein all its subgroups
admit non-constant holomorphic functions. Thus in our case
$\mathrm{M}^0=0$.
Moreover since $(\C^{\ast})^k$ and
$\mathrm{M}$ are Abelian their maximal tori are unique and it
follows that the maximal torus $(S^1)^s$ of $\mathrm{M}$ is
included into the maximal torus $(S^1)^k$ of $(\C^{\ast})^k$.
Therefore if $\dim_{\R} (\mathrm{M}\cap (S^1)^k)=0$ then
$\mathrm{M}$ is isomorphic to $\C^l$. Finally note that if an
injective morphism $\beta:\C^l\rightarrow (\C^{\ast})^k$ verifies
$\dim_{\R}(\IM (\beta)\cap (S^1)^k)=0$ then $l\leq k/2$.
This proves \textbf{(a)}. Part \textbf{(b)} is proved in a similar way.
\end{proof}

\begin{lem}\label{lem:u3bis}
Let \emph{M} be a complex subgroup of $(\C^{\ast})^{q}$
isomorphic to $\C^{p}$ and assume that we are in one of the
following two cases:
\begin{enumerate}[\bf (a)]\sep
\item $q=2k$, $p = k$ and $\dim_{\R} (\mathrm{M}\cap
(S^1)^{2k})=0$ or
\item $q=2k + 1$, $p = k+1$ and $\dim_{\R} (\mathrm{M}\cap
(S^1)^{2k+1})=1$.
\end{enumerate}
Then $\overline{\mathrm{M}}^{Zar}=(\C^{\ast})^{q}$, where
$\overline{\mathrm{M}}^{Zar}$ denotes the Zariski closure of
$\mathrm{M}$.
\end{lem}

\begin{proof}
Case \textbf{(a)} being similar and simpler we only give the proof
for case \textbf{(b)}. First notice that every connected algebraic
subgroup of $(\C^{\ast})^n$ is of the form $(\C^{\ast})^m$ (c.f.
\cite{Oni2}). Since $\overline{\mathrm{M}}^{Zar}$ is an algebraic
subgroup of $(\C^{\ast})^{2k+1}$ it is isomorphic to $(\C^{\ast})^l$
with $l\leq 2k+1$. Moreover its maximal torus $(S^1)^l$ is included
in the maximal torus $(S^1)^n$ of $(\C^{\ast})^n$. Now it is enough
to notice that $\mathrm{M}\cong\C^{k+1}$ cannot be immersed into
$(\C^{\ast})^l$ if $l<2k+1$ since this would contradict the
inequality
$$
\begin{array}{rcl}
2k+2+l & = & \dim_{\R}\C^{k+1} + \dim_{\R}(S^1)^l \\
& \leq & \dim_{\R}(\C^{\ast})^l+\dim_{\R}(\C^{k+1}\cap (S^1)^l)=2l+1.
\end{array}
$$
\end{proof}

\begin{prop}\label{lem:u4bis}
Let $\mathrm{K}$ be a compact connected semisimple Lie group with
Lie algebra $\mathfrak{k}$ endowed with a left-invariant complex
structure defined by a subalgebra $\mathfrak{l}$ of $\mathfrak{g}$
or with a nacs defined by a pair $\mathfrak{l}\subset \mathfrak{l}'$ of
subalgebras of $\mathfrak{g}$. Let $\mathrm B$ be a Borel subgroup of
$\mathrm G$ whose Lie algebra $\mathfrak{b}$ contains  $\mathfrak{l}$,
and $\mathfrak{l}'$ in the odd-dimensional case. Then $\mathfrak{t}:=
\mathfrak{b}\cap
\mathfrak{k}$ is a maximal Abelian subalgebra of $\mathfrak{k}$ and
the Cartan subalgebra $\mathfrak{r}=\mathfrak{t}^{\C}\subset
\mathfrak{b}$  fulfills:
\begin{enumerate}[\bf (i)]\sep
\item If $\dim_{\R}\mathfrak{k}=2n$ then
$\mathfrak{l}+\mathfrak{r}=\mathfrak{b}$ and
$\dim_{\R}(\mathfrak{l}\cap \mathfrak{r})=\rank \, \mathrm{K}$.
\item If $\dim_{\R}\mathfrak{k}=2n+1$ then
$\mathfrak{l}'+\mathfrak{r}=\mathfrak{b}$ and
$\dim_{\R}(\mathfrak{l}\cap \mathfrak{r})=\rank \, \mathrm{K}+1$.
\end{enumerate}
\end{prop}

\begin{proof}
Let us begin by proving that $\mathfrak{t}:=\mathfrak{b}\cap
\mathfrak{k}$ is a maximal Abelian subalgebra of $\mathfrak{k}$.
As the identity component of $\mathrm{B}\cap \mathrm{K}$ is a
closed connected Lie subgroup of a compact Lie group $\mathrm{K}$,
it is compact and it follows that $\mathfrak{b}\cap
\mathfrak{k}=[\mathfrak{b}\cap \mathfrak{k},\mathfrak{b}\cap
\mathfrak{k}]\oplus C(\mathfrak{b}\cap \mathfrak{k})$, where
$[\mathfrak{b}\cap \mathfrak{k},\mathfrak{b}\cap \mathfrak{k}]$ is
semisimple and $C(\mathfrak{b}\cap \mathfrak{k})$ is the center. On
the other hand, $\mathfrak{b}\cap \mathfrak{k}$ is solvable,
therefore it cannot admit a semisimple subalgebra and we conclude
that the subalgebra $\mathfrak{t}$ is Abelian.
Since all the Borel groups have the same dimension we have
$\dim_{\R}\mathfrak{b}=\dim_{\R}\mathfrak{k}+\rank \, \mathrm{K}$.
Therefore the inequality
$$\dim_{\R}\mathfrak{b}+\dim_{\R}\mathfrak{k}=
\dim_{\R}(\mathfrak{b}+\mathfrak{k})+\dim_{\R}(\mathfrak{b}\cap
\mathfrak{k})\leq
2\dim_{\R}\mathfrak{k}+\dim_{\R}(\mathfrak{b}\cap \mathfrak{k})$$
implies $\dim_{\R}(\mathfrak{b}\cap \mathfrak{k})\geq \rank \,
\mathrm{K}$.
So $\mathfrak{t} = \mathfrak{b}\cap \mathfrak{k}$ must be a maximal
Abelian subalgebra of $\mathfrak{k}$. Then
$\mathfrak{r}=\mathfrak{t}^{\C}$ is a Cartan subalgebra of
$\mathfrak{g}$ contained in $\mathfrak{b}$ and
$\dim_{\R}\mathfrak{r}=2\,\rank \, \mathrm{K}$.

Assume now that
$\dim_{\R}\mathfrak{k}=2n$ and $\rank \, \mathrm{K}=2r$. Since
$\dim_{\R}\mathfrak{l}=2n$ one has
$$
2n+2r =\dim_{\R}\mathfrak{b} \geq \dim_{\R}(\mathfrak{l} + \mathfrak{r})=2n+
4r-\dim_{\R}(\mathfrak{l}\cap \mathfrak{r}),
$$
so $2r\leq \dim_{\R}(\mathfrak{l}\cap \mathfrak{r})$. On the other
hand, since $\mathfrak{l}\cap \mathfrak{k}=0$ we have
$\dim_{\R}(\mathfrak{l}\cap \mathfrak{r})\leq 2r$ and this
concludes the proof in the even-dimensional case.

Finally let us consider the case
$\dim_{\R}\mathfrak{k}=2n+1$ and $\rank \, \mathrm{K}=2r+1$. Note
that $\langle \xi \rangle_{\R}=\mathfrak{l}'\cap
\mathfrak{k}\subset \mathfrak{b}\cap \mathfrak{k}=\mathfrak{t}$,
therefore $\mathfrak{l}'\cap \mathfrak{t}=\langle
\xi\rangle_{\R}$. By hypothesis
$\mathfrak{l}'+\mathfrak{r}\subset \mathfrak{b}$, and since
$\dim_{\R}\mathfrak{l}'=2n+2$ we have
$$
2n+2r+2 =\dim_{\R}\mathfrak{b}\geq \dim_{\R} (\mathfrak{l}'+\mathfrak{r})=2n+4+4r-\dim_{\R}(\mathfrak{l}'\cap\mathfrak{r}).
$$
So $2r+2\leq \dim_{\R}(\mathfrak{l}'\cap \mathfrak{r})$. Thus
$\mathfrak{l}'\cap \mathfrak{r}$ is a real subspace of
$\mathfrak{t}\oplus \mathrm{i}\mathfrak{t}$ of
dimension at least $2r+2$ whose intersection with $\mathfrak{t}$ has dimension
1, because $\mathfrak{l}'\cap \mathfrak{t}=\langle \xi \rangle_{\R}$.
Using Grassman formula one deduces that
$\dim_{\R}(\mathfrak{l}'\cap \mathfrak{r})=2r+2$. Now the
previous inequality implies
$\dim_{\R}(\mathfrak{l}'+\mathfrak{r})=2n+2r+2$ and consequently
$\mathfrak{l}'+\mathfrak{r}=\mathfrak{b}$.
\end{proof}

\begin{proof}(Theorem \ref{teo:borel})
We give the proof for the odd-dimensional case. That is,
we assume that $\dim_{\R}\mathrm{K}=2n+1$, $\rank \, \mathrm{K}=2r+1$
and that $\mathfrak{l}\subset \mathfrak{l}'$ are a pair of subalgebras
of $\mathfrak{g}$ defining a nacs on $\mathrm{K}$. For the
even-dimensional case the argument is analogous.

Proposition \ref{lem:u4bis} tells us that
$\mathfrak{t}=\mathfrak{b}\cap \mathfrak{k}$ is a maximal Abelian
subalgebra of $\mathfrak{k}$ and therefore
$\mathfrak{r}=\mathfrak{t}^{\C}$ is a Cartan subalgebra
contained in $\mathfrak{b}$. We must then show that
$\mathfrak{l}'=(\mathfrak{l}'\cap \mathfrak{r})\oplus
\mathfrak{u}$ and $\mathfrak{l}=(\mathfrak{l}\cap
\mathfrak{r})\oplus \mathfrak{u}$. Using the root decomposition of
$\mathfrak{g}$ by respect to $\mathfrak{r}$ we can write
$$\mathfrak{g}=\mathfrak{r}\oplus_{\alpha\in \Phi} \mathfrak{g}_{\alpha},
\quad \mathfrak{b}=\mathfrak{r}\oplus_{\alpha\in \widetilde{\Phi}}
\mathfrak{g}_{\alpha},\quad
\mathfrak{u}=\oplus_{\alpha\in{\widetilde{\Phi}}}\mathfrak{g}_{\alpha},$$
where $\Phi$ is the set of roots of $\mathfrak{g}$ relative to
$\mathfrak{r}$ different from zero, $\widetilde{\Phi}$ is a suitable
subset of $\Phi$ and
$\mathfrak{g}_{\alpha}$ are proper subspaces.

We consider now the
action on $\mathfrak{l}'$ of its Abelian subalgebra
$\mathfrak{l}'\cap \mathfrak{r}$. Since this algebra is Abelian
the endomorphisms of $\mathfrak{l}'$ defined
by the elements of $\mathfrak{l}'\cap \mathfrak{r}$ diagonalize
simultaneously. Thus we obtain a decomposition of $\mathfrak{l}'$
as a direct sum of eigenspaces
$$
\mathfrak{l}'=
\mathfrak{l}'_0 \oplus_{\alpha\in (\mathfrak{l}')^{\ast}-\{0\}} \mathfrak{l}'_{\alpha},
$$
where $\mathfrak{l}'_{\alpha}=\{X\in \mathfrak{l}': \,
[R,X]=\alpha(R)\cdot X, \;\; \forall R\in \mathfrak{l}'\cap
\mathfrak{r}\}$. Note that $\mathfrak{l}'_{\alpha}\subset
\mathfrak{u}$ for $\alpha\neq 0$. Indeed if there exists $R\in
\mathfrak{r}\cap \mathfrak{l}'$ such that $\alpha(R)\neq 0$
then for a given
$X\in \mathfrak{l}'_{\alpha}$ we have
$$
X= \frac{1}{\alpha(R)}[R,X]\in [\mathfrak{l}',\mathfrak{l}']
\subset [\mathfrak{b},\mathfrak{b}]=\mathfrak{u}.
$$
Moreover it is clear that $\mathfrak{l}'\cap \mathfrak{r}\subset
\mathfrak{l}'_0$. Now we want to prove that
$\mathfrak{l}'=(\mathfrak{l}'\cap \mathfrak{r})\oplus
\mathfrak{u}$ with $\mathfrak{l}_0'=\mathfrak{l}'\cap
\mathfrak{r}$ and $\oplus_{\alpha\in
(\mathfrak{l}')^{\ast}-\{0\}}\mathfrak{l}'_{\alpha}=\mathfrak{u}$.
In particular this will imply $\mathfrak{l}=(\mathfrak{l}\cap
\mathfrak{r})\oplus \mathfrak{u}$ ending the proof.
In fact, since
$[\mathfrak{l}',\mathfrak{l}']\subset \mathfrak{l}$ the same
argument above tells us
that $\mathfrak{l}'_{\alpha}\subset \mathfrak{l}$ for $\alpha\neq
0$ and therefore $\mathfrak{u}\subset \mathfrak{l}$. Then we deduce
$\mathfrak{l}=(\mathfrak{l}\cap \mathfrak{r})\oplus \mathfrak{u}$.

We will first see that $\mathfrak{l}'_0\subset
\mathfrak{r}$, which yields $\mathfrak{l}'\cap
\mathfrak{r}=\mathfrak{l}'_0$ and then we will conclude by an
argument of dimensions. We must check that given $X\in
\mathfrak{l}'$ such that $[R,X]=0$ for every $R\in
\mathfrak{r}\cap \mathfrak{l}'$ then $[R,X]=0$ for every $R\in
\mathfrak{r}$. Let $\mathrm{L}'$ and $\mathrm{H}$ be the connected
Lie subgroups of $\mathrm{G}$ corresponding to the Lie subalgebras
$\mathfrak{l}'$ and $\mathfrak{r}$ respectively and define $\mathrm{S}'$ as the
connected component of the identity in $\mathrm{L}'\cap \mathrm{H}$.
Recall that $\mathrm{H}\cong (\C^{\ast})^{2r+1}$.
According to  Lemma~\ref{lem:u2}  there are two
possibilities. In the first one $\mathrm{S}'\cong \C^{r}\times
\C^{\ast}$ and if we denote by $\mathrm{M}$ the subgroup of
$\mathrm{S}'$ isomorphic to $\C^{r}$ then $\dim_{\R}
\mathrm{M}\cap (S^1)^{2r+1}=0$.
In the second
one, $\mathrm{S}'\cong \C^{r+1}$ and $\dim_{\R}(\C^{r+1}\cap
(S^1)^{2r+1})=1$. In both cases lemma~\ref{lem:u3bis} implies
$\overline{\mathrm{S'}}^{Zar}=\mathrm{H}$. By hypothesis
$\mathrm{Ad}_h X=0$ for every $h\in \mathrm{L}'\cap
\mathrm{H}$ and $X\in \mathfrak{l}'$. Since this is an algebraic condition
and $\overline{\mathrm{S}'}^{Zar}=\mathrm{H}$ we deduce that
$\mathrm{Ad}_h X=0$ is also fulfilled for each $h\in \mathrm{H}$ and $X\in
\mathfrak{l}'$, and therefore
$\mathfrak{l'}_0\subset \mathfrak{r}$. We conclude that
$\mathfrak{l}'=(\mathfrak{l}'\cap\mathfrak{r})\oplus_{\alpha\in
\mathfrak{l'}^{\ast}-\{0\}}\mathfrak{l}'_{\alpha}$ with
$\mathfrak{l'}_{\alpha}\subset \mathfrak{u}$. Finally set
$\mathfrak{u}'=\oplus_{\alpha\in
(\mathfrak{l}')^{\ast}-\{0\}}\mathfrak{l}'_{\alpha}\subset \mathfrak{u}$
and notice that by Proposition \ref{lem:u4bis}
$$\dim_{\C}\mathfrak{u}'=n+1-r-1=n-r=\dim_{\C}\mathfrak{u}.$$
\end{proof}

The following proposition characterizes left-invariant CR-structures
on odd-di\-men\-sional Lie groups underlying
a nacs.

\begin{prop}\label{prop:CRsensecamp}
Let $\mathrm{K}^{2n+1}$ be an odd-dimensional semisimple compact
Lie group endowed with a left-invariant CR-structure of maximal dimension
defined by a complex subalgebra $\mathfrak{l}$ of $\mathfrak{g}$.
There exists a subalgebra $\mathfrak{l}'$ of $\mathfrak{g}$ such that the pair
$\mathfrak{l}\subset\mathfrak{l}'$ defines a nacs on $\mathrm{K}$
if and only if there is a Borel subalgebra $\mathfrak{b}$ of $\mathfrak{g}$
containing $\mathfrak{l}$ and such that
$\mathfrak{u}=[\mathfrak{b},\mathfrak{b}]$ is contained in $\mathfrak{l}$.
In that case the choices can be made in such a way that
$\mathfrak{l'}\subset\mathfrak{b}$.
\end{prop}

\begin{proof}
The ``if" part of the statement is just Theorem \ref{teo:borel}.
Assume that there is a Borel subalgebra $\mathfrak{b}$
containing $\mathfrak{l}$, which already implies that $\mathfrak{l}$ is solvable,
and such that  $\mathfrak{u}=[\mathfrak{b},\mathfrak{b}]
\subset \mathfrak{l}$. One can show as in Proposition~\ref{lem:u4bis}
that if we set $\mathfrak{t} = \mathfrak{b} \cap \mathfrak{k}$ then
$\mathfrak{r}=\mathfrak{t}^{\C}\subset\mathfrak{b}$ is a Cartan subalgebra
of $\mathfrak{g}$.
Moreover one has $\mathfrak{b}=\mathfrak{r}\oplus \mathfrak{u}$. The hypothesis
$\mathfrak{u}\subset \mathfrak{l}$ then implies
$\mathfrak{l}=(\mathfrak{l}\cap \mathfrak{r})\oplus \mathfrak{u}$. Since $\dim_{\C}\mathfrak{b}=n+r+1$
and $\dim_{\C}\mathfrak{l}=n$ the last two equalities imply
$\dim_{\C}(\mathfrak{l}\cap\mathfrak{r})<\dim_{\C}\mathfrak{r}$.
Therefore there exists $\xi\in \mathfrak{r}\backslash (\mathfrak{l}\cap
\mathfrak{r})$. We define $\mathfrak{l}'=\mathfrak{l}\oplus
\langle \xi \rangle$ and we claim that the pair $\mathfrak{l}\subset
\mathfrak{l}'$ defines a nacs. Clearly $\dim_{\R}(\mathfrak{l}'\cap
\mathfrak{k})=1$ and moreover $\mathfrak{l}$ is an ideal of
$\mathfrak{l}'$. Indeed each $X\in \mathfrak{l}$ can be written
$X=X_0+X_1$ where $X_0\in \mathfrak{l}\cap \mathfrak{r}$ and $X_1\in
\mathfrak{u}$, thus $[\xi, X]=[\xi, X_0]+[\xi,X_1] = [\xi,X_1] \in
\mathfrak{u}\subset \mathfrak{l}$.
\end{proof}


\section{A geometric construction of invariant structures}

In this section we introduce a geometrical construction which provides
left-invariant complex structures or nacs on a given semisimple compact
Lie group $\mathrm{K}$. It is shown that each invariant structure on
$\mathrm{K}$ is obtained in that way.  As a corollary we also prove
the existence of left-invariant complex structures or nacs on any compact
Lie group. This was already known by the work of Samelson and Wang
(cf. \cite{Samel} and \cite{Wang2}) in the case of complex structures
and the work of  Charbonel and Khalgui on CR-structures (cf. \cite{Charb})
in the case of nacs. The interest of the geometrical approach is that
after a slight modification it can also produce non-invariant, and therefore
different, structures on $\mathrm{K}$. The generalization of this construction
will be carried out in Section~4.

Let us assume that the compact Lie group $\mathrm{K}$ is semisimple
and let us fix a maximal torus $\mathrm{T}\cong (S^1)^{\rank \, \mathrm{K}}$
of $\mathrm{K}$ and a Borel
subgroup $\mathrm{B}$ of $\mathrm{G}$ containing $\mathrm{T}$.
Let $\mathfrak{t}$ and $\mathfrak{b}$ be the Lie algebras of  $\mathrm{T}$
and $\mathrm{B}$ respectively. The connected Lie subgroup $\mathrm{H}$ of $\mathrm{G}$ associated to the Lie algebra $\mathfrak{r}=\mathfrak{t}^{\C}$
is a Cartan subgroup and it is isomorphic to $(\C^{\ast})^{\rank \, \mathrm{K}}$.
One has $\mathrm{B}=\mathrm{H}\cdot \mathrm{U}$ where $\mathrm{U}$
is the nilpotent Lie group associated to
$\mathfrak{u}=[\mathfrak{b},\mathfrak{b}]$. We also recall that the Iwasawa decomposition states $\mathrm{G}\cong\mathrm{K}\cdot \mathrm{A}
\cdot \mathrm{U}$
where $\mathrm{A}\cong\R^{\rank \, \mathrm{K}}$ is the connected Lie subgroup associated to the Lie
subalgebra $\mathrm{i} \mathfrak{t}$.

A Lie group morphism $\Lambda:\C^{l}\rightarrow
(\C^{\ast})^{q}$ is the composition with the exponential map of
a $\C$-linear morphism
$$
\Lambda^0:\C^l\rightarrow \C^q
$$
$$z=(z_1,\ldots,z_l)^t\mapsto \mathrm{M}\cdot z$$
where $\mathrm{M}=(m_i^j)$ is a $(q \times l)$-complex matrix. Let
$\textbf{e}_1,...,\textbf{e}_{2l}$ be the canonical basis of
$\R^{2l}\cong \C^l$, i.e. $\textbf{e}_1=(1,0,\dots,0)$,
$\textbf{e}_2=(\mathrm{i},0,\dots,0)$, \dots, $\textbf{e}_{2l}=(0,\dots,\mathrm{i})$,
and denote by $\mathrm{A}_{\Lambda}$ the $(q\times 2l)$-real matrix
whose colums are $\RE(\Lambda^0(\textbf{e}_i))$. That is
$$
\mathrm{A}_{\Lambda}=\left(\begin{array}{ccccc}
\RE m_1^1 & -\IM m_1^1 & \ldots & \RE m_1^{l} & -\IM m_1^{l} \\
\vdots & \vdots & & \vdots & \vdots\\
\RE m_{q}^1 & -\IM m_{q}^1 & \ldots & \RE m_{q}^{l} &
-\IM m_{q}^{l}
\end{array} \right).
$$
We denote by $\mathrm {B}_{\Lambda}$ the real matrix of dimension
$q\times (2l-2)$ which is obtained from $\mathrm{A}_{\Lambda}$ just
by removing the first two columns and keeping the last $2l-2$ ones.

For us $q$ will be the rank of $\mathrm K$ and we consider two cases.
In the first one the dimension of $\mathrm K$ is even and $\rank \, \mathrm{K}=2r$
and $l=r$. In this case $\mathrm{A}_{\Lambda}$ is a square
$(2r\times 2r)$-real matrix and we have:

\begin{lem}\label{lem:cons1}
Let $\Lambda:\C^{r}\rightarrow \mathrm{H}$ be a Lie group morphism.
The Lie subgroup $\Lambda(\C^r)$ is
transverse to $\mathrm{T}$ if and only if
\begin{equation}
\mathrm{det} \mathrm{A}_{\Lambda}\, \neq 0. \tag{I}
\end{equation}
In that case $\Lambda(\C^r)$ is a closed subgroup of $\mathrm{H}$.
\end{lem}

The proof of this statement is similar and simpler than the one of
Lemma~\ref{lem:cons2} below and we omit it.

In the second case the dimension of $\mathrm K$ is odd and $\rank \, \mathrm{K}=2r+1$
and $l=r+1$. In this case the real matrices $\mathrm{A}_{\Lambda}$ and
$\mathrm{B}_{\Lambda}$ have respective dimensions $(2r+1)\times (2r+2)$
and $(2r+1)\times 2r$ and we have

\begin{lem}\label{lem:cons2}
Let $\Lambda:\C^{r+1}\rightarrow \mathrm{H}$ be a Lie group
morphism. Then the Lie subgroups $\Lambda(\{0\}\times\C^r)$ and $\Lambda(\C^{r+1})$ of $\mathrm{H}$ fulfill the transversality conditions
$$
\Lambda(\{0\}\times\C^r)\cap \mathrm{T}=\{0\},\quad \dim_{\R}\Lambda(\C^{r+1})\cap \mathrm{T}=1
$$
if and only if the matrices $\mathrm{A}_{\Lambda}$ and
$\mathrm{B}_{\Lambda}$ have maximal rank, i.e.
\begin{equation}
\mathrm{rank}\,\mathrm{A}_{\Lambda}=2r+1 \quad \mathrm{and} \quad \mathrm{rank}\,\mathrm{B}_{\Lambda}=2r \tag{II}
\end{equation}
In that case $\Lambda(\{0\}\times\C^r)$ is a closed subgroup of
$\mathrm{H}$.
\end{lem}

\begin{rem}
Under the hypothesis of the lemma the morphism $\Lambda$ is
injective when restricted to $\{0\}\times\C^r$ and $\dim \ker
\Lambda=0$.
\end{rem}

\begin{proof}
Notice that $\mathrm{T}$ is included in $\mathrm{H}$ as
$(S^1)^{2r+1}\subset (\C^{\ast})^{2r+1}$. We have
$$
\Lambda(z_1,...,z_{r+1})=   \left(\exp\bigg({\sum_{j=1}^{r+1} m_1^jz_j}\bigg),...,
\exp\bigg({\sum_{j=1}^{r+1} m_{2r+1}^j z_j}\bigg)\right)
$$
and therefore
$\Lambda(0,z_2,...,z_{r+1})$ intersects $\mathrm{T}$
 if and only if
\begin{equation}\label{equat}
\RE\bigg(\sum_{j=2}^{r+1} m_1^jz_j\bigg)=...=\RE\bigg(\sum_{j=2}^{r+1} m_{2r+1}^jz_j\bigg)=0.
\end{equation}
Setting $z_j=x_j+\mathrm{i}y_j$ condition (\ref{equat}) can be
rewritten as
$$
\sum_{j=2}^{r+1} \bigg(\RE m_k^j x_j-\IM m_k^jy_j\bigg)=0 \qquad \forall k=1,...,2r+1.
$$
The coefficients of this homogeneous system are the
entries of the matrix $\mathrm{B}_{\Lambda}$. Therefore system
(\ref{equat}) admits a unique solution $z_2=...=z_{r+1}=0$ if and only if
$\mathrm{rank}\,\mathrm{B}_{\Lambda}=2r$. Similarly
the condition $\dim_{\R}\Lambda(\C^{r+1})\cap
\mathrm{T}=1$ is equivalent to $\mathrm{rank}\,
\mathrm{A}_{\Lambda}=2r+1$.

Finally we prove that in this situation
$\Lambda(\{0\}\times \C^r)$ is closed in
$\mathrm{H}=(\C^{\ast})^{2r+1}$. It is enough to check that it is
closed in a neighborhood of the identity $e = \Lambda(0)$.
We are assuming that the rank of $\mathrm{B}_{\Lambda}$ is maximal
and therefore
$\Lambda(z)$ is close to $e$ if and only if
$$
\bigg|\sum_{j=2}^{r+1} \bigg(\RE m_k^j x_j-\IM m_k^jy_j\bigg)\bigg|<\epsilon \qquad
\forall k=1,...,2r+1.
$$
Since the restriction to $\{0\}\times \C^r$ of the morphism $\Lambda$
is an immersion this proves that
$\Lambda(\{0\}\times \C^r)$ is closed in $\mathrm{H}$.
\end{proof}

We assume that $\Lambda:\C^{l}\rightarrow \mathrm{H}$, where  $l=r$ or $l=r+1$
according to the parity of $\dim_{\R}\mathrm{K}$, is a Lie group morphism fulfilling the
transversality conditions (I) or (II). Since $\mathrm{B} = \mathrm{H}\cdot \mathrm{U}$
is the normalizer $N(\mathrm{U})$ of $\mathrm{U}$ we can conclude the
following. If $\dim_{\R}\mathrm{K}= 2n$ and therefore its rank is also even,
$\rank\,\mathrm{K}= 2r$, then $\mathrm{L}_{\Lambda}=
\Lambda(\C^{r})\cdot \mathrm{U}$ is a well defined
complex Lie subgroup of $\mathrm{G}$ of complex dimension $n$. We denote by
$\mathrm{l}_{\Lambda}$ the corresponding Lie subalgebra.
In a similar way, if $\dim_{\R}\mathrm{K}= 2n+1$ and
$\rank\,\mathrm{K}= 2r+1$ then $\mathrm{L}_{\Lambda}=
\Lambda(\{0\}\times \C^{r})\cdot \mathrm{U}$ and $\mathrm{L}'_{\Lambda}=
\Lambda(\C^{r+1})\cdot \mathrm{U}$
are Lie subgroups of $\mathrm{G}$ of complex dimensions $n$ and
$n+1$ respectively. We denote by
$\mathrm{l}_{\Lambda}$ and $\mathrm{l}'_{\Lambda}$ the corresponding
Lie subalgebras.

Under the transversality hypothesis (I), or (II), see Lemmas \ref{lem:cons1},
and \ref{lem:cons2} respectively,
the Lie subgroup  $\mathrm{L}_{\Lambda}$
is closed in $\mathrm{G}$ although this will not be necessarily the case for
$\mathrm{L}'_{\Lambda}$. This is proved in next theorem which states that
each such morphism $\Lambda$ induces a left-invariant complex structure or
nacs on $\mathrm{K}$ and that every such an invariant structure is obtained in this way.

\begin{thm}\label{teo:cons2}
Let $\mathrm{K}$ be a semisimple compact connected Lie group, and let
$\mathrm{H}$ be a Cartan subgroup of $\mathrm{G}$ and $\mathrm{B}$ a
Borel subgroup containing $\mathrm{H}$.
Assume that $\Lambda:\C^{l}\rightarrow \mathrm{H}$, where $l=r$ or
$l=r+1$ according to the parity of $\dim_{\R}\mathrm{K}$, is a Lie group
morphism fulfilling the transversality
condition $\mathrm{(I)}$ or $\mathrm{(II)}$.  Then the group $\mathrm{L}_{\Lambda}$ is closed in
$\mathrm{G}$ and the Lie subalgebra $\mathrm{l}_{\Lambda}$, or the pair
$\mathrm{l}_{\Lambda}\subset \mathrm{l}'_{\Lambda}$, define a left-invariant
complex structure, or nacs, on $\mathrm{K}$. More precisely
\begin{enumerate}[\bf (i)]\sep
\item
if $\dim_{\R}\mathrm{K}= 2n$ and $\rank\,\mathrm{K}= 2r$ then the natural inclusion
$\mathrm{K}\hookrightarrow \mathrm{G}/\mathrm{L}_{\Lambda}$ is a biholomorphism.
\item
if $\dim_{\R}\mathrm{K}= 2n+1$ and $\rank\,\mathrm{K}= 2r+1$ then
$\mathrm{K}\hookrightarrow \mathrm{G}/\mathrm{L}_{\Lambda}$ is a CR-embedding.
Moreover there is holomorphic vector field $\zeta$ on
$\mathrm{G}/\mathrm{L}_{\Lambda}$ transverse to $\mathrm{K}$ and whose real part
is tangent to $\mathrm{K}$ and defines the CR-action of the nacs on $\mathrm{K}$.
\end{enumerate}
Conversely, every left-invariant complex structure or nacs on $\mathrm{K}$
is induced by such a morphism $\Lambda:\C^{l}\rightarrow \mathrm{H}$ for a
suitable choice of $\mathrm{T}$ and $\mathrm{B}$.
\end{thm}

The above construction of invariant complex structures or nacs on $\mathrm{K}$
can also be described in the following way. The inclusion $\mathrm{K}\subset\mathrm{G}$
induces an embedding of $\mathrm{K}$
into the complex manifold $\mathrm{G}/\mathrm{U}$ that, using representation
theory, can be shown to be a Zariski open subset of
an affine algebraic variety.
On the other hand, since $\mathrm{H}$
is contained in the normalizer $N(\mathrm{U}) = \mathrm{B}$ of $\mathrm{U}$ the
morphism $\Lambda:\C^{l}\rightarrow \mathrm{H}$ induces a well defined action of
$\C^{l}$ on the quotient $\mathrm{G}/\mathrm{U}$. The key point is that conditions
(I) and (II) imply that this action is transverse to $\mathrm{K}$. Therefore we can apply
the following basic fact

\begin{lem}\label{lem:folcomplexa}
Let ${\mathcal F}$ be a holomorphic foliation on a complex manifold
$\mathrm{X}$. A real submanifold $\mathrm{M}$ of $\mathrm{X}$ transverse to
${\mathcal F}$ inherits a transversely holomorphic foliation in a natural way.
If $\dim\mathrm{M} = {\mathrm{cod}_{\R}\mathcal{F}}$ then the transversely
holomorphic foliation is just a complex structure on $\mathrm{M}$.
\end{lem}

In our situation $\mathrm{X}=\mathrm{G}/\mathrm{U}$ and the foliation ${\mathcal F}$ is defined by
the orbits of the $\C^{l}$-action. There are two cases. If $\dim_{\R}\mathrm{K}= 2n$ then
$l = r$ and the dimension of $\mathrm{K}$ is precisely the real codimension of ${\mathcal F}$.
Therefore ${\mathcal F}$ induces a complex structure on $\mathrm{K}$. In fact $\mathrm{K}$ is
naturally identified to the orbit space of the action. The second case, when $\dim_{\R}\mathrm{K}= 2n+1$
and $l=r+1$, is more involved. In this case ${\mathcal F}$ induces on $\mathrm{K}$ a transversely holomorphic foliation ${\mathcal{F}_{\mathrm{K}}}$ of dimension one. The CR-structure on $\mathrm{K}$ is determined by $\Lambda(\{0\}\times \C^{r})$ and there is a fundamental
vector field of the
$\C^{r+1}$-action
whose real flow $\varphi_s$ is tangent to $\mathrm{K}$, preserves the CR-structure and defines
the foliation ${\mathcal{F}_{\mathrm{K}}}$.

\begin{rem}
Since all Borel subgroups of $\mathrm{G}$ are conjugated,
different choices of the maximal torus $\mathrm{T}$ and the
Borel group $\mathrm{B}$ in the above theorem produce
left-invariant complex structures, or nacs, on $\mathrm{K}$ which are conjugated.
That is, up to conjugacy, the only way
of obtaining different left-invariant structures is by
making different choices of the morphism $\Lambda:
\C^{l}\rightarrow \mathrm{H}$
\end{rem}

\begin{proof}(of the Theorem)
We give the proof in the odd-dimensional case. Notice that
$\mathrm{B}$, $\mathrm{H}$, and $\mathrm{U}$ are closed in $\mathrm{G}$
and that the map $\varphi:\mathrm{H}\times \mathrm{U}\rightarrow \mathrm{B}=
\mathrm{H}\cdot\mathrm{U}$, given by $\varphi(h,u)=h\cdot u$,
is a diffeomorphism. Therefore if $\{h_m\}$, $\{u_m\}$ are sequences in
$\Lambda(\{0\}\times \C^r)$ and $\mathrm{U}$ respectively with
$h_m\cdot u_n\rightarrow b\in \mathrm{B}$ then we can write  $b=h\cdot u$ with
$h\in \mathrm{H}$ and $u\in \mathrm{U}$ in a unique way and one has
$h_m\rightarrow h$, $u_m\rightarrow u$. By Lemma \ref{lem:cons2}
the group $\Lambda(\{0\}\times \C^r)$ is closed in $\mathrm{H}$ and hence
$h\in \Lambda(\{0\}\times \C^r)$ proving that $\mathrm{L}_{\Lambda}$ is closed in
$\mathrm{G}$.

The map $\mathrm{K}\rightarrow \mathrm{G}/\mathrm{L}_{\Lambda}$
induced by the inclusion $\mathrm{K}\subset \mathrm{G}$ is an embedding
of K as a real hypersurface of the complex manifold
$\mathrm{G}/\mathrm{L}_{\Lambda}$. This
inclusion induces a CR-structure on K which is left-invariant and which is
determined by the Lie subalgebra $\mathfrak{l}_{\Lambda}$ of $\mathfrak{g}$.
Notice also that
$\dim_{\C}\mathfrak{l}_{\Lambda}=n$,
$\dim_{\C}\mathfrak{l}_{\Lambda}'=n+1$,  $\mathfrak{l}_{\Lambda}\cap
\mathfrak{k}=0$, $\dim_{\R}\mathfrak{l}_{\Lambda}'\cap
\mathfrak{k}=1$ and that
$[\mathfrak{l}_{\Lambda}',\mathfrak{l}_{\Lambda}']\subset
[\mathfrak{b},\mathfrak{b}]=\mathfrak{u}\subset
\mathfrak{l}_{\Lambda}$, which implies that the pair of subalgebras
$\mathfrak{l}_{\Lambda}\subset\mathfrak{l}'_{\Lambda}$ define a
left-invariant nacs on $\mathrm{K}$.
Moreover the vector field $\zeta$ on $\mathrm{G}/\mathrm{L}_{\Lambda}$
induced by a non-zero element $\xi\in\mathfrak{l}_{\Lambda}'\cap
\mathfrak{k}$ has the required properties. (Up to a linear coordinate change in
$\C^{r+1}$ one can assume that $\xi = \Lambda(\textbf{e}_1)$.)

The converse follows from Theorem \ref{teo:borel} and the observation
that the only pairs $\mathrm{L}\subset \mathrm{L}'$ of complex Lie
subgroups of $\mathrm{B}$ containing  $\mathrm{U}$ and such that
$\mathrm{T}\cap
\mathrm{L}=\{e\}$ and $\mathrm{T}\cap \mathrm{L}'=\langle \xi
\rangle_{\R}$ are those of the form
$\mathrm{L}_{\Lambda}=\Lambda(\{0\}\times \C^r)\cdot
\mathrm{U}\subset \mathrm{L'}_{\Lambda}=\Lambda'(\C^{r+1})\cdot
\mathrm{U}$ where $\Lambda:\C^{r+1}\rightarrow (\C^{\ast})^{2r+1}$
is a Lie group morphism verifying the transversality condition (II).
\end{proof}

\begin{exam}\label{exam:su2}
Let us consider the semisimple Lie group $\mathrm{K}=\mathrm{SU}(2)$
whose universal complexification is $\mathrm{SL(2,\C)}$. Its Lie algebra can
be written $\mathfrak{k} = \mathfrak{su}(2) = \langle e_1, e_2, e_3\rangle_{\R}$
with $[e_1, e_2]=2e_3$, $[e_2, e_3]=2e_1$ and $[e_3, e_1]=2e_2$. We can take
as Cartan group $\mathrm{H}$ and Borel group $\mathrm{B}$ the Lie subgroups associated to the subalgebras
$$
\mathfrak{h}=\langle e_1 \rangle_{\C}\quad \mathfrak{b}=\langle e_1,e_2+\mathrm{i}e_3 \rangle_{\C}.
$$
Then $\mathfrak{u}=\langle u_2+\mathrm{i} u_3 \rangle_{\C}$. Since
$\dim_{\C}\mathfrak{h}=1$, up to conjugation there is only one left-invariant
nacs on $\mathrm{SU}(2)$, the one defined by the pair of subalgebras
$\mathfrak{l}= \mathfrak{u}$ and $\mathfrak{l}'= \mathfrak{b}$.

Nevertheless using Proposition \ref{prop:CRsensecamp}  one sees
that there are left-invariant
CR-structures on $\mathrm{SU}(2)$ which do not underlie a nacs.
One can take for instance any algebra of the family
$\mathfrak{l}_\alpha=\langle e_1+
\alpha (e_2+\mathrm{i}e_3)\rangle_{\C}$ for $\alpha \in \R^*$.
The same considererations can be made for SO(3) as $\mathfrak{so}(3,\R)\cong
\mathfrak{su}(2)$.
\end{exam}

We end this section by considering the general case of a compact Lie group
$\mathrm K$ which is not necessarily semisimple. As a corollary of
Theorem \ref{teo:cons2} we prove the well known result (cf. \cite{Samel},  \cite{Wang2} and \cite{Charb}) of the existence of a left invariant complex structure or nacs on $\mathrm K$. Before proving this we make some easy considerations.

\begin{rem}\label{rem:nacsgeneral}
Let $\mathrm K$ be a semisimple compact connected Lie group of
even dimension endowed with a left-invariant complex structure
defined by a subalgebra $\mathfrak{l}$ of $\mathfrak{g}$.  Then
this invariant structure induces a left-invariant nacs on the product
$\mathrm K\times S^1$ defined by the pair
$\mathfrak{l}\subset \mathfrak{l}'=\mathfrak{l}\oplus
\langle \derpp{t}\rangle^{\C}$ where $\langle \derpp{t}\rangle^{\C}$
is the vector field determined by the natural $S^1$-action.

In a similar way one can see that if ${\mathrm K}_1$ and ${\mathrm K}_2$
are compact connected Lie groups carrying respectively a left-invariant nacs
and a left-invariant complex structure then the product
${\mathrm K}_1 \times {\mathrm K}_2$ is endowed with a natural left-invariant
nacs.
\end{rem}

\begin{cor}\label{cor:geninv}
Let $\mathrm{K}$ be a compact connected Lie group. Then $\mathrm{K}$
admits left-invariant complex structures, or nacs, according to the parity of
its dimension.
\end{cor}

\begin{proof}
Let us consider only the odd-dimensional case. Assume first that
$\mathrm{K}$ is semisimple. Since there exist group morphisms
$\Lambda:\C^{r+1}\rightarrow \mathrm{H}$
fulfilling the transversality condition (II), Theorem~\ref{teo:cons2}
shows the existence of left-invariant nacs on $\mathrm{K}$.

Let us consider now the general case. As it was explained in
\ref{rem:generalcase} the compact Lie group
$\mathrm{K}$ can be written as a quotient $\mathrm{K} =
\Gamma\backslash\mathrm{K}'\times (S^1)^p$ where $\mathrm{K}'$
is semisimple and $\Gamma$ is a finite subgroup of
$\mathrm{K}'\times (S^1)^p$. It follows from remarks \ref{rem:nacsxs1}
and \ref{rem:nacsgeneral} that there also exists left-invariant nacs on
the product $\mathrm{K}'\times (S^1)^p$. Since $\Gamma$ is contained
in the center they induce invariant nacs on the quotient $\mathrm{K}$.
\end{proof}


\section{Construction of non-invariant structures}

The aim of this section is to generalize the construction made in Section~3 in order to
define non-invariant complex structures or nacs on the semisimple compact
Lie group $\mathrm{K}$. As it is explained after the statement of
Theorem~\ref{teo:cons2} that construction can be understood as follows. There is
 a natural embedding of $\mathrm{K}$ into the complex manifold
 $\mathrm{G}/\mathrm{U}$. On the other hand  $\mathrm{G}/\mathrm{U}$ can be
 endowed with a holomorphic foliation $\mathcal F$ whose leaves are the orbits
 of a suitable $\C^l$-action. This action is associated to the choice of an Abelian
 subgroup $\Lambda(\C^l)$ of $\mathrm{H}\subset\mathrm{G}$ and, if it fulfills
 a certain transversality condition (conditions (I) or (II)), then it induces either a left-invariant
complex structure or a nacs on $\mathrm{K}$.

In this section we construct more general foliations $\mathcal F$ on
$\mathrm{G}/\mathrm{U}$
also induced by $\C^l$-actions but not necessarily associated to an Abelian
subgroup of $\mathrm{G}$. Under certain conditions $\mathcal F$
induces a complex structure or a nacs on $\mathrm{K}$ which in general it is not
invariant. In the even-dimensional case this condition is simple.  Due to
Lemma~\ref{lem:folcomplexa} transversality of $\mathcal F$ with $\mathrm{K}$
already assures the existence of a complex structure. In the odd-dimensional case
the condition is necessarily stronger because transversality assures the existence of a transversely holomorphic flow on $\mathrm{K}$ but it is not sufficient by itself to determine neither the CR-structure on $\mathrm{K}$ neither to guarantee the existence of a transverse $\R$-action.

\medskip

As above we fix a Cartan subgroup $\mathrm{H}$ of $\mathrm{G}$, a Borel
subgroup $\mathrm{B}$ containing $\mathrm{H}$ and we denote by $\mathrm{U}$
the derived group of $\mathrm{B}$, i.e. $\mathrm{U}= [\mathrm{B},\mathrm{B}]$.
With these choices the Iwasawa decomposition can be written
$\mathrm{G}=\mathrm{K}\cdot\mathrm{A}\cdot\mathrm{U}$ where $\mathrm{A}$
is a sugroup of $\mathrm{H}$ isomorphic to $\R^{\rank\,\mathrm{K}}$.
Therefore there is a fibration
$$
\begin{CD}
\mathrm{K} & @>>> & \mathrm{G}/\mathrm{U} \\
&&&&@VV\pi V \\
&&&& \mathrm{A}
\end{CD}
$$
whose fibers are the translations $\mathrm{K} a$ of $\mathrm{K}$ by the elements of
$\mathrm{A}$ (since $\mathrm{A}\subset N(\mathrm{U})$ this notation is not ambiguous).
We put    $\mathrm{K} = \mathrm{K} e$.

We consider first the even-dimensional case, that is $\dim_{\R}\mathrm{K}= 2n$ and $\rank\,\mathrm{K}= 2r$. Given a locally free holomorphic $\C^{r}$-action
$\varphi:\C^{r}\times \mathrm{G/U}\rightarrow \mathrm{G/U}$ we will denote by
$F_x$ the orbit of $[x]\in \mathrm{G/U}$, that is $F_x=\varphi(\C^{r},[x])$.

\begin{defn}
We say that a locally free holomorphic $\C^{r}$-action
$\varphi:\C^{r}\times \mathrm{G/U}\rightarrow \mathrm{G/U}$ fulfills the transversality
condition (III$_a$), where $a\in \mathrm{A}$,  if the orbits of the action are transverse
to the fiber $\mathrm{K} a$, that is
\begin{equation}
\dim_{\R}(F_p\cap \mathrm{K} a)=0 \quad\text{for each }  p\in \mathrm{K} a. \tag{$\mathrm{III}_a$}
\end{equation}
In case $a=e$ we put $\mathrm{(III_e)}=\mathrm{(III)}$.
\end{defn}

We consider now the odd-dimensional case. In this situation
$\dim_{\R}\mathrm{K}= 2n+1$ and $\rank\,\mathrm{K}= 2r+1$. We recall that
a CR-structure with a transverse $\R$-action (or nacs) on $\mathrm{K}$
is determined by a real vector field $\xi$ on $\mathrm{K}$ and two
complex subbundles $V$ and $V'$ of $T^{\C}\mathrm{K}$ of respective
ranks $n$ and $n+1$ such that: (i) $V_p\cap T_p\mathrm{K}=0$ and
 $V'=V\oplus \langle \xi \rangle_{\C}$, (ii) $V$ and $V'$ are involutive, and
 (iii) $[\xi,V]\subset V$.

Given a locally free holomorphic
$\C^{r+1}$-action $\varphi:\C^{r+1}\times \mathrm{G/U}\rightarrow
\mathrm{G/U}$ and $x\in \mathrm{G}$ we denote by $F_x$ and $F'_x$
the orbits of $[x]\in \mathrm{G/U}$ by the action of the groups $\{0\}\times \C^{r}$
and $\C^{r+1}$ respectively. That is
$$
F_x=\varphi(\{0\}\times \C^{r},[x]) \qquad F'_x=\varphi(\C^{r+1},[x]).
$$
We denote by $z_0,z_1,\ldots,z_r$ the linear coordinates of
$\C^{r+1}$.

\begin{defn}\label{trans4}
A locally free holomorphic $\C^{r+1}$-action $\varphi:\C^{r+1}\times
\mathrm{G/U}\rightarrow \mathrm{G/U}$ is said to fulfill the transversality
condition $\mathrm{(IV_a)}$, where $a\in \mathrm{A}$,  if there exist $\lambda\in\C$
such that for each $p\in \mathrm{K} a$ one has
\begin{equation}
\left\{
\begin{array}{ll}
\rm{(i)} & \dim_{\R}(F_p\cap \mathrm{K} a)=0  \text{ and } \dim_{\R}(F'_p\cap \mathrm{K} a)=1 \\
\rm{(ii)} & \xi=d\varphi\big(\RE\big(\lambda\derpp{z_0}\big) \big)\text{ is tangent to }
\mathrm{K} a.
\end{array}
\right.
\tag{$\mathrm{IV}_a$}
\end{equation}
In case $a=e$ we put $\mathrm{(IV_e)}=\mathrm{(IV)}$.
\end{defn}

The following statement is straightforward.

\begin{lem}\label{lem:folCR}
A locally free holomorphic $\C^{r+1}$-action $\varphi:\C^{r+1}\times
\mathrm{G/U}\rightarrow \mathrm{G/U}$ fulfilling the transversality condition
$\mathrm{(IV)}$ induces a nacs on $\mathrm{K}$.
\end{lem}

The next result, which is also straightforward, is the key for constructing more
general $\C^l$-actions on $\mathrm{G/U}$.

\begin{lem}
The action of $(\C^{\ast})^{2l}\cong \mathrm{H}\times
\mathrm{H}=\mathrm{B}/\mathrm{U}\times \mathrm{B}/\mathrm{U}$ on the
homogeneous space $\mathrm{G}/\mathrm{U}$ given by
$(h_1,h_2,[g])\mapsto [h_1 \cdot g\cdot h_2]$
is well defined.
\end{lem}

Let $\Lambda:\C^{l}\rightarrow \mathrm{H\times H}$ be a given
Lie group morphism. We write $\Lambda=(\Lambda_1,\Lambda_2)$ where
$\Lambda_i:\C^{l}\rightarrow \mathrm{H}$ are the components of $\Lambda$.
The composition of such a morphism with the above $(\mathrm{H\times H})$-action
induces a $\C^{l}$-action on $\mathrm{G/U}$ that we denote
$$
\varphi_{\Lambda}:\C^{l}\times \mathrm{G/U}\rightarrow\mathrm{G/U}.
$$

\begin{rem}
It is easy to see that the action $\varphi_{\Lambda}$ associated to a morphism
of the form $\Lambda=(e,\Lambda_2)$ fulfills condition (III) or (IV) if and only if $\Lambda_2$
fulfills condition (I) or (II) respectively.
\end{rem}

As above we set $l=r$ if $\rank \, \mathrm{K} = 2r$ and $l = r+1$ if
$\rank \, \mathrm{K} = 2r + 1$. With this notation we have

\begin{prop}\label{prop:pictureCR}
Let $\Lambda:\C^{l}\rightarrow\mathrm{H\times H}$ be a
Lie group morphism inducing a locally free action
$\varphi_{\Lambda}:\C^{l}\times \mathrm{G/U}\rightarrow \mathrm{G/U}$
which fulfills the transversality condition $\mathrm{(III)}$, or $\mathrm{(IV)}$. Then
$\varphi_{\Lambda}$ also fulfills $\mathrm{(III_a)}$, or $\mathrm{(IV_a)}$,
for each $a\in \mathrm{A}$.
\end{prop}

\begin{proof}
Since $\mathrm{A}\subset N(\mathrm{U})$ one has
\begin{equation}\label{transit}
\varphi_{\Lambda}(c,[x\, a])=[\Lambda_1(c)\, x\, a \,\Lambda_2(c)]=
[\Lambda_1(c)\, x\,  \Lambda_2(c)\, a]=
\varphi_{\Lambda}(c,[x])\, a.
\end{equation}
Let us consider the odd-dimensional case. Given $[y]\in\mathrm{K} a$
there are $x\in \mathrm{K}$ and
$v\in \mathrm{U}$ such that $y=x\, a \, v$. Now conditions (i) in
Definition~\ref{trans4} follow directly from  (\ref{transit})
and condition (ii) follows from the differential version of (\ref{transit}).
\end{proof}

\begin{thm}\label{prop:CRnoinv}
Let $\mathrm{K}$ be a semisimple compact connected Lie group.
Every morphism of Lie groups
$\Lambda:\C^{l}\rightarrow \mathrm{H}\times\mathrm{H}$
inducing a locally free holomorphic action
$\varphi_{\Lambda}:\C^{l}\times \mathrm{G/U}\rightarrow
\mathrm{G/U}$ and fulfilling the transversality condition
$\mathrm{(III)}$, or $\mathrm{(IV)}$, determines on
$\mathrm{K}$ a complex structure, or a nacs, in a natural way. Moreover,
such a structure is left-invariant if and only if
$\Lambda=(e,\Lambda_2)$.
\end{thm}

The fact that a morphism $\varphi_{\Lambda}$ fulfilling $\mathrm{(III)}$ or
$\mathrm{(IV)}$ induces on $\mathrm{K}$ a complex structure or a nacs follows from Lemmas \ref{lem:folcomplexa} and \ref{lem:folCR}. The characterization of
the invariant structures is given by the next proposition which is only stated
for the odd-dimensional case. The analogous statement for the even-dimensional
case can be proved in a similar way.

\begin{prop}\label{prop:criinvCR}
Let $\Lambda:\C^{r+1}\rightarrow \mathrm{H\times H}$ be a Lie group
morphism inducing a locally free holomorphic action
$\varphi_{\Lambda}:\C^{r+1}\times \mathrm{G/U}\rightarrow
\mathrm{G/U}$ and fulfilling the transversality condition $\mathrm{(IV)}$.
Then the following conditions are equivalent:
\begin{enumerate}[\bf (a)]\sep
\item The nacs on \emph{K} induced by $\Lambda$ is left-invariant.
\item One has $g\cdot F_y=F_{g\cdot y}$
and $g\cdot F_y'=F_{g\cdot y}'$ for each $y,g\in \mathrm{G}$.
\item The morphism $\Lambda$ is of the form $\Lambda=(e,\Lambda_2)$.
\end{enumerate}
\end{prop}


\begin{proof}
Let us begin by proving $\bf{(a)}\Rightarrow \bf{(b)}$. Assume that the
nacs on $\mathrm{K}$ is left-invariant. In particular the two subbundles
$V$ and $V'=V\oplus \langle \xi \rangle_{\C}$ of $T^{\C}\mathrm{K}$ defining
the structure are invariant by $\mathrm{K}$. This means that
\begin{equation}\label{distribution}
\begin{array}{lll}
k\cdot \mathrm{d}\varphi_{\Lambda}(\{0\}\times\C^r,T_{[x]} \mathrm{G/U}) &=&
\mathrm{d}\varphi_{\Lambda}(\{0\}\times\C^r,T_{[k\, x]} \mathrm{G/U}),\\
k\cdot \mathrm{d}\varphi_{\Lambda}(\C^{r+1},T_{[x]} \mathrm{G/U}) &=&
\mathrm{d}\varphi_{\Lambda}(\{0\}\times\C^{r+1},T_{[k\, x]} \mathrm{G/U}).
\end{array}
\end{equation}
for every $x,k\in \mathrm{K}$ (by abuse of notation we also denote by
$k\cdot$ the differential map of left multiplication by $k$). Denote by
$\mathcal F$ and $\mathcal F'$ the foliations on $ \mathrm{G/U}$
whose leaves are $F_y$ and $F'_y$ respectively and let and $\mathcal D$
and $\mathcal D'$ be the corresponding tangent distributions.
Then (\ref{distribution})
says that the restriction,  $\mathcal D|_\mathrm{K}$
and $\mathcal D' |_\mathrm{K}$, of these distributions to $ \mathrm{K}$ are
invariant by the action of $\mathrm{K}$. Now the differential
version of the identity (\ref{transit}) imply that the whole distributions
$\mathcal D$ and $\mathcal D'$ are in fact invariant by $\mathrm{K}$.
By integration we deduce that the foliations themselves are
$\mathrm{K}$-invariant, that is
\begin{equation}\label{leaves}
k\cdot F_y=F_{k\, y}\quad \text{ and } \quad k\cdot F_y'=F_{k\, y}'
\end{equation}
for each $y\in \mathrm{G}$ and $k\in \mathrm{K}$.

Now we denote by $\widetilde{\mathrm{G}}$ the set of elements in
$\mathrm{G}$ preserving the holomorphic foliations $\mathcal F$ and
$\mathcal F'$. That is
$$
\widetilde{\mathrm{G}}=\{g\in \mathrm{G}\,|\,g\cdot F_y=F_{g\,
y},\; g\cdot F_y'=F_{g\, y}',\; \forall y\in \mathrm{G}\}.
$$
Note that $\widetilde{\mathrm{G}}$ is a complex Lie subgroup of $\mathrm{G}$
because it is defined by holomorphic conditions and one has
$e\in \widetilde{\mathrm{G}}$ and $g\cdot F_{g^{-1}\,
y}=F_y$, so $F_{g^{-1}\cdot y}= g^{-1}\cdot F_y$. Since $\mathrm{G}$ is totally
real in $\mathrm{G}$ we conclude that
$\mathrm{G}=\widetilde{\mathrm{G}}$ which proves $\bf{(b)}$.
\medskip

It is clear that condition $\bf{(c)}$ implies $\bf{(a)}$ therefore it is
sufficient to prove the implication $\bf{(b)}\Rightarrow \bf{(c)}$.

Let
$c\in \C^{r+1}$ be given. We want to see that $\Lambda_1(c)= e$.  Denote by
$\mathfrak{r}$, $\mathfrak{b}$ and $\mathfrak{u}$ the Lie subalgebras
associated to $\mathrm{H}$,  $\mathrm{B}$ and $\mathrm{U}$
respectively. Using the Cartan decomposition we can write, as above,
$\mathfrak{b}=\mathfrak{r} \oplus_{\alpha\in \widetilde{\Phi}}
\mathfrak{g}_{\alpha}$ and $\mathfrak{u}=\oplus_{\alpha\in
\widetilde{\Phi}}\mathfrak{h}_{\alpha}$. Then
$$
\mathfrak{b}'=\mathfrak{r} \oplus_{\alpha\in
\widetilde{\Phi}} \mathfrak{g}_{-\alpha}
$$
is a Lie subalgebra whose corresponding connected Lie group is
also a Borel group  $\mathrm{B}'$. If $ \mathrm{U}'$ is its derived group
then $\mathrm{B}'=\mathrm{H\cdot U'}$ and
$\mathrm{U}\cap \mathrm{U}'=\{e\}$.

Let us choose $g\in \mathrm{U}'$. By
hypothesis there exists $d\in \C^{r+1}$ such that
$$
g^{-1}\cdot \varphi_{\Lambda}(c,[g])=\varphi_{\Lambda}(d,[e]).
$$
Equivalently there exists $u\in \mathrm{U}$ such that
\begin{equation}\label{diamondsuit}
g\,\Lambda_1(d)\, \Lambda_2(d)=
\Lambda_1(c)\, g \, \Lambda_2(c)\, u .
\end{equation}
Note that $g\in \mathrm{U}'$ and
$\Lambda_1(c),\Lambda_2(c),\Lambda_1(d),\Lambda_2(d)\in \mathrm{H}$,
therefore $g\,\Lambda_1(d)\, \Lambda_2(d)$ and
$\Lambda_1(c)\, g \, \Lambda_2(c)$ belong to $\mathrm{B}'$. As
$\mathrm{U}\cap \mathrm{B}'=\{e\}$ equation (\ref{diamondsuit}) implies
\begin{equation}\label{heartsuit}
g\, \Lambda_1(d)\,\Lambda_2(d)=\Lambda_1(c)\, g \, \Lambda_2(c)
\end{equation}
which can be rewritten as
$$
\Lambda_1(d)^{-1}\, g\, \Lambda_1(d)=\Lambda_1(c-d)\, \Lambda_2(c-d)
\, \Lambda_2^{-1}(c-d)\, g \, \Lambda_2(c-d)
$$
where $\Lambda_1(d)^{-1}\, g\,
\Lambda_1(d),\Lambda_2^{-1}(c-d)\, g \, \Lambda_2(c-d)\in
\mathrm{U}'$ and $\Lambda_1(c-d)\, \Lambda_2(c-d)\in \mathrm{H}$.
Since $\mathrm{H}\cap \mathrm{U}'=\{e\}$ this implies
$\Lambda_1(c-d)=\Lambda_2(d-c)$ or equivalently
$$
\Lambda_1(c)\, \Lambda_2(c)=\Lambda_1(d)\, \Lambda_2(d).
$$
Combining the last equation with (\ref{heartsuit}) we obtain
$$
g=\Lambda_1(c)^{-1}\, g \, \Lambda_1(c).
$$
This last identity is verified for each $c\in \C^{r+1}$ and each $g\in\mathrm{U}'$,
and thus for each $g\in\mathrm{B}'= \mathrm{H}\cdot\mathrm{U}'$.
Therefore it implies that the image of
$\Lambda_1$ belongs to the center of $\mathrm{B}'$
which is equal to the center $Z(\mathrm{G})$ of $\mathrm{G}$
(cf. \cite{Hump2}, p.140). Since  $Z(\mathrm{G})=Z(\mathrm{K})$ and
$\mathrm{K}$ is semisimple the center is finite.
By continuity $\Lambda_1(c) =\Lambda_1(0)=e$ as desired.
\end{proof}

\begin{rem}\label{rem:center}
It follows from its definition that the action $\varphi_{\Lambda}$ of ${\mathbb C}^l$
on $\mathrm{G/U}$ is $Z(\mathrm{G})$-equivariant. This means that
one has
$g\cdot F_y=F_{g\,y}$ and $g\cdot F_y'=F_{g\, y}'$ for each $y\in\mathrm{G}$
and each $g\in Z(\mathrm{G})=Z(\mathrm{K})$. So, although the complex structures
or nacs given by Theorem~\ref{prop:CRnoinv} are not
$\mathrm{K}$-invariant in general, all of them are invariant by the center $Z(\mathrm{K})$
of $\mathrm{K}$.

\end{rem}

We end the section deducing from Theorem~\ref{prop:CRnoinv}
the existence of non-invariant complex structures or nacs on each
compact Lie group $\mathrm{K}$ different from $S^1$. First we discuss the
particular cases $\mathrm{K}=\mathrm{SU}(2)$ and $\mathrm{SO}(3,\R)$.

\begin{exam}
In this example we describe the nacs on
$\mathrm{K} =\mathrm{SU}(2)\cong S^3$ that are obtained
by means of the construction given in the above theorem.
In particular we show the existence of non-invariant nacs.
With the choices made in Example~\ref{exam:su2},
the Iwasawa decomposition of $\mathrm{G}=\mathrm{SL}(2,\C)$ is
$$
\mathrm{SL}(2,\C)=\mathrm{SU}(2)\cdot \mathrm{A}\cdot \mathrm{U}
$$
where $\mathrm{A}=\left\{\left( \begin{smallmatrix} \lambda & 0 \\ 0
& \lambda^{-1} \end{smallmatrix}\right): \, \lambda\in \R^{+}
\right\}$, $\mathrm{H}=\left\{\left(
\begin{smallmatrix} \alpha & 0 \\ 0 & \alpha^{-1}
\end{smallmatrix}\right): \, \alpha\in \C^{\ast} \right\}$,
$\mathrm{U}=\left\{\left( \begin{smallmatrix} 1 & a \\ 0 & 1
\end{smallmatrix}\right): \, a\in \C \right\}$.
The linear action of $\mathrm{SL}(2,\C)$ on $\C^2$ identifies
$\mathrm{G/U} = \mathrm{SL}(2,\C)/\mathrm{U}\cong \C^2\backslash \{0\}$.
With this identification $\mathrm{SU}(2)$ is embedded as the unit sphere
in $ \C^2\backslash \{0\}$. Then the action of
$\mathrm{H}\times\mathrm{H}\cong (\C^{\ast})^2$ on
$\mathrm{G/U}\cong\C^2\backslash \{0\}$ corresponds to
$$
(\mathrm{H}\times\mathrm{H})\times
\C^2\backslash \{0\}\rightarrow \C^2\backslash \{0\}
$$
$$(\alpha,\beta),(z,w)\mapsto (\alpha\beta z,\alpha^{-1}\beta w).$$
A morphism $\Lambda:\C\rightarrow (\C^{\ast})^2$ is of the form
$\Lambda(t)=(e^{at},e^{bt})$ for $a,b\in \C$ and the induced
$\C$-action $\varphi_{\Lambda}$ on $\C^2\backslash \{0\}$ is
$$\C \times \C^2\backslash \{0\}\rightarrow \C^2\backslash \{0\}$$
$$(t,(z,w))\mapsto (e^{(a+b)t}z,e^{(b-a)t}w).$$
Notice that $\varphi_{\Lambda}$ is induced by the vector field
$$\eta=(a+b)z\derpp{z}+(b-a)w\derpp{w},$$ which intersects $S^3$
in a 1-dimensional orbit if
and only if there does not exist $\mu\in \R^{-}$ such that
$a+b=\mu(b-a)$.
In fact a straightforward computation shows that the action
$\varphi_{\Lambda}$ fulfills the transversality condition $\mathrm{(IV)}$
if and only if there exists $\mu\in \R^{+}$ such that $a+b=\mu(b-a)$.
In each of these cases the action defines a nacs on
$\mathrm{SU}(2)$. A direct computation shows also that the structure is
left-invariant if and only if $a=0$.

We remark that all these structures are invariant by the involution
$\nu(z,w) = (-z,-w)$ and therefore induce nacs on $\mathrm{SO}(3,\R)$,
which is the quotient of $\mathrm{SU}(2)$ by its center.
\end{exam}

\begin{cor}\label{teo:gennoninv}
Let $\mathrm{K}$ be a compact connected Lie group with $\dim\mathrm{K}>1$.
Then $\mathrm{K}$ admits non-invariant complex structures or nacs
(and in particular a non-invariant CR-structure of maximal
dimension).
\end{cor}

\begin{proof}
Assume first that the group $\mathrm{K}$ is semisimple. If
$\rank\,\mathrm{K}>1$, one can easily construct non-invariant structures as
small deformations of invariant ones. For instance, in the odd-dimensional
case, it is enough to determine the action of $\Lambda$ on the vector
$\textbf{e}_1=(1,0,\dots,0)$ in such a way that condition (\textbf{ii}) in
the transversality hypothesis (IV) is verified and then to chose the images of the
other vectors in the basis in such a way that the CR-structure is not invariant.

On the other hand the only compact
connected Lie groups of rank 1 are $S^1$, $\mathrm{SO(3)}$ and
$\mathrm{SU(2)}$. The group $S^1$ only admits an invariant nacs
but it is excluded by the hypothesis. For the other groups the existence of
non-invariant nacs has been shown before.

In the general case the statement can be proved in a way similar
to Corollary~\ref{cor:geninv} taking into account Remarks~\ref{rem:generalcase}
and \ref{rem:center}.
\end{proof}


\bibliographystyle{amsplain}


\providecommand{\bysame}{\leavevmode\hbox to3em{\hrulefill}\thinspace}
\providecommand{\MR}{\relax\ifhmode\unskip\space\fi MR }
\providecommand{\MRhref}[2]{%
  \href{http://www.ams.org/mathscinet-getitem?mr=#1}{#2}
}
\providecommand{\href}[2]{#2}

\end{document}